\documentclass{amsart}
\usepackage{amsfonts}
\usepackage{amscd}
\usepackage{amsmath}

\usepackage[all]{xy}
\CompileMatrices

\theoremstyle{plain} \numberwithin{equation}{section}
\newtheorem{theorem}{Theorem}[section]

\newtheorem{lemma}[theorem]{Lemma}

\theoremstyle{definition}

\newtheorem{remark}[theorem]{Remark}

\newtheorem{question}{Question}
\numberwithin{equation}{section}

\newcommand{\so}{\mathfrak{so}}
\newcommand{\su}{\mathfrak{su}}
\newcommand{\ad}{\mathrm{ad}}

\newcommand{\fa}{\mathfrak{a}}

\newcommand{\fc}{\mathfrak{c}}
\newcommand{\fg}{\mathfrak{g}}
\newcommand{\fh}{\mathfrak{h}}
\newcommand{\fk}{\mathfrak{k}}
\newcommand{\fl}{\mathfrak{l}}
\newcommand{\fm}{\mathfrak{m}}
\newcommand{\fn}{\mathfrak{n}}
\newcommand{\fp}{\mathfrak{p}}

\newcommand{\fr}{\mathfrak{r}}

\newcommand{\ft}{\mathfrak{t}}

\newcommand{\fz}{\mathfrak{z}}

\newcommand{\R}{\mathbb{R}}

\newcommand{\N}{\mathbb{N}}

\newcommand{\Ad}{\mathop{\mathrm{Ad}} }

\newcommand{\oline}{\overline}

\newcommand{\rank}{\mathop{\rm rank}}

\begin{document}
\title[Iwasawa decompsition]{A novel characterization of the Iwasawa decomposition of a simple 
Lie group}

\author{Bernhard Kr\"otz}
\thanks{MPI, DFG}
\address{TU Darmstadt, MPI}


\maketitle
This appendix is about (essential) uniqueness of the 
{\it Iwasawa (or horospherical)  decomposition} $G=KAN$of 
a  semisimple Lie group $G$. This means: 

\begin{theorem}\label{t=A} Assume that $G$ is a connected Lie group
with simple Lie algebra $\fg$. Assume that $G=KL$ for some 
closed subgroups $K, L<G$ with $K\cap L$ discrete. Then up to order, 
the Lie algebra $\fk$ of $K$ is maximally compact, and the Lie algebra
$\fl$ of $L$ is isomorphic to $\fa+\fn$, the Lie algebra 
of $AN$. 
\end{theorem}

\section{General facts on decompositions of Lie groups}
\label{s-one}

For a group $G$, a subgroup $H<G$ and an element $g\in G$ we 
define $H^g=gHg^{-1}$. 

\begin{lemma}\label{l-1} Let $G$ be a group and $H,L<G$ subgroups. 
Then the following statements are equivalent: 
\begin{enumerate} \item $G=HL$ and $H\cap L=\{ {\bf 1}\}$. 
 \item $G=HL^g$  and $H\cap L^g=\{ {\bf 1}\}$ for all $g\in G$.
\end{enumerate}
\end{lemma}

\begin{proof} Clearly, we only have to show that $(ii)\Rightarrow (i)$. 
Suppose that $G=HL$ with $H\cap L=\{ {\bf 1}\}$. 
Then we can write $g\in G$ as $g=hl$ for some $h\in H$ and 
$l\in L$. Observe that $L^g=L^h$ and so 
$$H\cap L^g= H\cap L^h=H^h\cap L^h=(H\cap L)^h=\{{\bf 1}\}\ .$$
Moreover we record 
$$HL^g=HL^h=HLh=Gh=G\, .$$
\end{proof}

In the sequel, capital Latin letters will denote 
real Lie groups and the corresponding lower case fractur
letters will denote the associated Lie algebra, i.e. $G$ is 
a Lie group with Lie algebra $\fg$.

\begin{lemma} \label{l-2}Let $G$ be a Lie group and $H,L<G$ closed 
subgroups. Then the following statements are 
equivalent: 
\begin{enumerate} \item $G=HL$ with $H\cap L=\{{\bf 1}\}$. 
\item The multiplication map 
$$H\times L \to G, \qquad (h,l)\mapsto hl$$ 
is an analytic diffeomorphism.
\end{enumerate}
\end{lemma}
\begin {proof} Standard structure theory. \end{proof}
 
If $G$ is a Lie group with closed subgroups $H,L<G$ such that 
$G=HL$ with $H\cap L=\{{\bf 1}\}$, then we refer to $(G,H,L)$
as a {\it decomposition triple}.

\begin{lemma}\label{l-3} Let $(G, H,L)$ be a decomposition 
triple.  Then: 
\begin{equation}\label{eq=ac}
(\forall g\in G) \qquad \fg=\fh +\Ad(g)\fl \qquad \hbox{and} \qquad 
\fh\cap \Ad(g)\fl=\{0\}\, .\end{equation} 
\end{lemma}

\begin{proof} In view of Lemma \ref{l-2}, the map 
$H\times L \to G, \ \ (h,l)\mapsto hl$ is 
a diffeomorphism. In particular, the differential at 
$({\bf 1}, {\bf 1})$ is a diffeomorphism which means 
that $\fg=\fh+\fl$, $\fh\cap\fl=\{0\}$. As we may replace 
$L$ by $L^g$, e.g. Lemma \ref {l-1}, the assertion follows.
\end{proof}

\begin{question} Assume that $G$ is connected. 
Is it then true that $(G,H,L)$ is a decomposition 
triple if and only if the algebraic condition (\ref{eq=ac}). 
is satisfied.
\end{question}

\begin{remark} If the Lie algebra $\fg$ splits
into a direct sum of subalgebras $\fg=\fh+\fl$, then 
we cannot conclude in general that $G=HL$ holds. 
For example, let $\fh= \fm +\fa +\fn$ is a minimal 
parabolic subalgebra and $\fl =\oline {\fn}$ is the opposite of $\fn$.  
Then $HL=MAN\oline N$ is the open Bruhat cell in $G$. 
A similar example is when $\fg=\mathrm{sl}(n,\R)$ 
with $\fh$ the upper triangular matrices and 
$\l =\so(p,n-p)$ for $0<p<n$. In this case $HL\subset G$ is 
a proper open subset. Notice that in both examples 
condition (\ref{eq=ac}) is violated as $\fh\cap \Ad(g)\fl
\neq \{0\}$ for appropriate $g\in G$. 
\end{remark} 
 
\section{The case of one factor being maximal compact}

Throughout this section $G$ denotes a semi-simple connected 
Lie group with associated Cartan decomposition
$\fg=\fk+\fp$. Set $K=\exp \fk$ and note that 
$\Ad (K)$ is maximal compact subgroup in $\Ad(G)$. 

\par For what follows we have to recall some results 
of Mostow  on maximal solvable subalgebras in $\fg$. 
Let $\fc\subset \fg$ be a Cartan subalgebra. 
Replacing $\fc$ by an appropriate $\Ad(G)$-conjugate
we may assume that $\fc=\ft_0 +\fa_0$ with 
$\ft_0\subset \fk$ and $\fa_0\subset \fp$. Write 
$\Sigma=\Sigma(\fa, \fg)\subset \fa^*\backslash\{0\}$ for the 
non-zero $\ad \fa_0$-spectrum on $\fg$.  For $\alpha\in \Sigma$
write $\fg^\alpha$ for the associated eigenspcae.
Call 
$X\in \fa_0$ {\it regular} if $\alpha(X)\neq 0$ for all 
$\alpha\in \Sigma$. Associated to a regular element $X\in \fa$
we associate a nilpotent subalgebra 
$$\fn_X=\bigoplus_{\alpha\in \Sigma\atop \alpha(X)>0} \fg^\alpha\, .$$ 
If $\fa\subset \fp$ happens to be maximal abelian, then we will write 
$\fn$ instead of $\fn_X$.

With this notation we have:

\begin{theorem} \label{th=1}Let $\fg$ be a semi-simple Lie algebra. Then 
the following assertions hold: 
\begin{enumerate}
\item Every maximal solvable subalgebra $\fr$ of $\fg$ contains
a Cartan subalgebra $\fc$ of $\fg$. 
\item Up to conjugation with an element of $\Ad(G)$ every 
maximal solvable subalgebra of $\fg$ is of the form 
$$\fr=\fc +\fn_X$$
for some regular element $X\in \fa_0$.
\end{enumerate} 
\end{theorem}
\begin{proof} \cite{M2}, Theorem 4.1.
\end{proof}

We choose a maximal abelian subspace $\fa\subset \fp$ 
and write $\Sigma=\Sigma(\fg, \fa)$ for the associated 
root system. For a choice of positive roots we obtain 
a unipotent subalgebra $\fn$. Write $\fm=\fz_\fk(\fa)$ and 
fix a Cartan subalgebra $\ft\subset \fm$.  
Write $A, N, T$  for the analytic subgroups of $G$ corresponding 
to $\fa$, $\fn$, $\ft$. Notice that $\ft+\fa+\fn$ is 
a maximal solvable subalgebra by Theorem \ref{th=1}.

\begin{lemma} Let $L<G$ be a closed subgroup such that $G=KL$ with 
$K\cap L=\{ {\bf 1}\}$. Then there is an Iwasawa decomposition 
$G=NAK$ such that 
\begin{equation} \label{comp} 
N\subset L \subset TAN \qquad \hbox{and}\qquad L\simeq LT/T\simeq AN\, .
\end{equation}
Conversely, if $L$ is a closed subgroup of $G$ satisfying (\ref{comp}), 
then $G=KL$ with $K\cap L=\{{\bf 1}\}$. 
\end{lemma}
\begin{proof} Our first claim is that $L$ contains 
no non-trivial compact subgroups. In fact, let $L_K \subset L$ 
be a compact subgroup. 
As all maximal compact subgroups of $G$ are conjugate, we find a $g\in G$ 
such that $L_K^g\subset K$. But $L_K^g\cap K \subset L^g\cap K =\{{\bf 1}\}$
by Lemma \ref{l-1}. This establishes our claim. 

\par Next we show that $L$ is solvable. For that let $L=S_L\times R_L$ be a Levi decomposition, where 
$S$ is semi-simple and $R$ is reductive. If $S\neq {\bf 1}$, then 
there is a non-trivial maximal compact subgroup $S_K\subset S$.
Hence $S={\bf 1}$ by our previous claim and $L=R_L$ is solvable.

\par Next we turn to the specific structure of $\fl$, the Lie algebra 
of $\fl$.  
Let $\fr=\fc+\fn_X$ be a maximal solvable subalgebra
of $\fg$ which contains $\fl$. As before we write 
${\fc}={\ft}_0+{\fa}_0$ for the Cartan subalgebra of $\fr$. . We claim 
that ${\fa}_0=\fa$ is maximal abelian in $\fp$. In fact, notice that 
$\fl\cap \ft_0=\{0\}$ and so 
$\fl \hookrightarrow\fr/\ft_0\simeq \fa_0 +\fn_X$ injects as vector spaces. 
Hence 
$$\dim \fl =\dim \fa +\dim \fn \leq \dim \fa_0 +\dim \fn_X. $$
But $\dim \fa_0\leq \dim \fa$ and $\dim \fn_X\leq \dim \fn$ 
and therefore $\fa=\fa_0$. Hence $\fr=\ft+\fa+\fn$. 
As $\fl\simeq \fr/\ft$ as vector spaces we thus get hat 
$L\simeq LT/T\simeq R/T\simeq AN$ as homogeneous spaces.  
We now show that $N\subset L$ which will follow from 
$\fn \subset [\fl,\fl]$. For that choose a regular element 
$X\in \fa$. By what we know already, we then find an element 
$Y\in \ft$ such that $X+Y\in \fl$. Notice that 
$\ad (X+Y)$ is invertible on $\fn$ and hence 
$\fn\subset [X+Y, \fn]$. Finally, observe that 
$$[X+Y, \fn]=[X+Y, \fr]=[X+Y,\fl+\ft]=[X+Y, \fl]$$
which concludes the proof of the first assertion of the 
lemma. 

\par Finally, the second assertion of the lemma is immediate 
from the Iwasawa decomposition of $G$. 
\end{proof}

\section{Manifold decompositions for decomposition triples}

Throughout this section $G$ denotes a connected Lie group. 
\par

\par Let $(G,H,L)$ be a decomposition triple and let us 
fix maximal compact subgroups 
$K_H$ and $K_L$ of $H$ and 
$L$ respectively.  We choose a maximal compact 
subgroup $K$ of $G$ such that $K_H\subset K$. 
As we are free to replace $L$ by any conjugate $L^g$, we 
may assume in addition that $K_L\subset K$.

We then have the following fact, see also \cite{O}, Lemma 1.2.

\begin{lemma} Let $(G,H,L)$ be a decomposition 
triple. Then $(K, K_H, K_L)$ is a decomposition triple, i.e.
the map 
$$K_H\times K_L\to K, \ \ (h,l)\mapsto hl$$
is a diffeomorphism.
\end{lemma}

Before we proof the Lemma we recall a fundamental 
result of Mostow 
concernig the topology of a connected Lie group $G$, cf. \cite{M1} 
If $K<G$ is a maximal compact subgroup of $G$, then 
there exists a a vector space $V$ and 
a homeomorphism $G\simeq K\times V$. In particular $G$ is
a deformation retract of $K$ and thus 
$H_\bullet  (G, \R) =H_\bullet (K,\R)$. 

\begin{proof} As $H\cap L=\{{\bf 1}\}$, it follows that 
$K_H\cap K_L=\{ {\bf 1}\}$. Thus compactness of 
$K_L$ and $K_H$ implies that 
the map 

$$K_H\times K_L\to K, \ \ (h,l)\mapsto hl$$
has closed image. It remains to show that the image is 
open. This will follow from $\dim K_H +\dim K_L=\dim K$. 
In fact $G\simeq H\times L$ implies that $G$ is 
homeomorphic to $K_H\times K_L \times V_H\times V_L$ for 
vector spaces $V_H$ and $V_L$. Thus  
$$H_\bullet (K, \R)=H_\bullet (G, \R)=H_\bullet (K_H\times K_L, \R)$$
and K\"unneth implies for any $n\in \N_0$ that 

$$H_n (K, \R)\simeq \sum_{j=0}^n H_j  (K_H,\R)\otimes H_{n-j}
(K_L, \R).$$ 
Now, for an orientable connected compact manifold $M$ we recall that 
$H_{\dim M} (M, \R)=\R$ and  $H_n(M,\R)=\{0\}$ for 
$n>\dim M$.  Next Lie groups are orientable and 
we deduce from the K\"unneth identity from above 
that   $\dim K_H + \dim K_L=\dim K$. This 
concludes the proof of the lemma. 
\end{proof}

Let us write $\fk_\fh$ and $\fk_\fl$ for the Lie algebras 
of $K_H$ and $K_L$ respectively. Then, as $(K, K_H, K_L)$ 
is a decomposition triple, it follows from Lemma \ref{l-3} that 
$$\fk=\fk_\fh+\Ad(k)\fk_\fl\qquad \hbox{and}
\qquad \fh\cap \Ad(k)\fl=\{0\}\, .$$ 
Let now $\ft_h\subset \fk_\fh$ be a maximal toral subalgebra and extend 
it to a maximal torus $\ft$, i.e. $\ft_\fh\subset \ft$. 
Now pick a maximal toral subalgebra $\ft_\fl$. Replacing 
$\fl$ by an appropriate $\Ad(K)$-conjugate, we may assume that 
$\ft_\fl\subset \ft$ (all maximal toral subalgebras in $\fk$ 
are conjugate). Finally write $T, T_H, T_L$ for the corresponding 
tori in $T$. 

\begin{lemma}\label{l-4} If $(K, K_H, K_L)$ is a decomposition 
triple for a compact Lie group $K$, then 
$(T, T_H, T_L)$ is a decomposition triple for the maximal 
torus $T$. 
In particular 
\begin{equation} \rank K=\rank K_H +\rank K_L\, .
\end{equation} 
\end{lemma}
 
\begin{proof} We already know that $\ft_\fh +\ft_\fl\subset 
\ft$ with $\ft_\fh\cap \ft_\fl=\{0\}$. It remain s to verify that 
$\ft_\fh+\ft_\fl=\ft$. 
We argue by contradiction. Let $X\in \ft, X\not\in \ft_\fh+\fh_\fl$. 
As $\fk=\fk_\fh+\fk_\fl$, we can write 
$X=X_\fh +X_\fl$ for some $X_\fh\in \fk_\fh$ and $X_\fl\in \fk_\fl$. 
\par For a compact Lie algebra $\fk$ with maximal toral subalgebra
$\ft\subset \fk$ we recall the direct vector space decomposition 
$\fk=\ft\oplus [\ft, \fk]$. As $\ft_\fh +\ft_\fl\subset 
\ft$ we hence may assume that $X_\fh \in [\ft_\fh, \fk_\fh]$ and 
$X_\fl \in [\ft_\fl, \fk_\fl]$.  
But then we get 
$$X=X_\fh+X_\fl\in  [\ft_\fh, \fk_\fh]+  [\ft_\fl, \fk_\fl]\subset 
[\ft, \fk]$$
and therefore $X\in \ft\cap [\ft, \fk]=\{0\}$, a contradiction. 
\end{proof}

\section {Decompositions of compact Lie groups}

Decompositions of compact Lie groups is an algebraic 
feature as the following Lemma, essentially due to Oni\v s\v cik, 
shows.

\begin{lemma} Let $\fk$ be a compact Lie algebra
and $\fk_1, \fk_2< \fk$ be two subalgebras. 
Then the following statements are equivalent: 
\begin{enumerate}
\item $\fk=\fk_1 +\fk_2$ with $\fk_1\cap\fk_2=\{0\}$
\item Let $K , K _1, K _2$ be simply connected 
Lie groups with Lie algebras $\fk, \fk_1$ and $\fk_2$. 
Write $\iota_i: K _i\to K $, $i=1,2$ for the natural 
homomorphisms sitting over the inclusions
$\fk_i\hookrightarrow \fk$. Then the map 
$$m: K _1\times K _2\to K ,
 \ \ (k_1, k_2)\mapsto \iota_1(k_1) \iota_2(k_2)$$
is a homeomorphism. 
\end{enumerate}
\end{lemma}

\begin {proof} The implication 
$(ii)\Rightarrow (i)$ is clear.
We establish $(i)\Rightarrow (ii)$. 
We need that $m$ is onto and deduce this from \cite{O}, Th. 3.1. 
Then $K$ becomes a homogeneous space for the left-right
action of $K_1\times K_2$. The stabilizer of $\bf 1$ is given
by the discrete subgroup $F=\{ (k_1, k_2): \iota_1(k_1)=\iota_2(k_2)^{-1}$, 
i.e. $K\simeq K_1\times K_2/F$. As 
$K_1$ and $K_2$ are simply connected, we conclude that 
$\pi_1(K) =F $ and thus $F=\{\bf 1\}$ as $K$ is simply 
connected.  
\end{proof}

We now show the main result of this 
section. 

\begin{lemma}\label{l-5} Let $(K, K_1, K_2)$ be a decomposition 
triple of a connected compact simple Lie group. Then 
$K_1={\bf 1}$ or $K_2={\bf 1}$. 
\end{lemma}

Before we prove this, a few remarks are in order. 

\begin{remark} (a) If $K$ is of exceptional type, then the result 
can be easily deduced from $\dim K=\dim K_1 +\dim K_2$ and the 
rank equality $\rank K =\rank K_1 +\rank K_2$, cf. Lemma \ref{l-4}. 
For example if $K$ is of type $G_2$. Then a non-trivial 
decomposition $K=K_1K_2$ must have $\rank K_i=1$, i.e.
$\fk_i=\su(2)$. But 
$$14=\dim K\neq \dim K_1 +\dim K_2=6\, .$$ 
\par (b) The assertion of the lemma is not true if we 
only require $K=K_1 K_2$ and drop $K_1\cap K_2=\{{\bf 1}\}$. 
For example if $K$ is of type $G_2$. then $K=K_1 K_2$ with 
$K_i$ locally $\mathrm{SU}(3)$ and $K_1\cap K_2 =T$ a maximal torus. 
\end{remark}

\begin{proof} The proof is short, but uses 
a powerful tool, namely the structure of the cohomology ring 
of the compact group $K$. See for instance \cite{Oz} or 
\cite{Kosz}. 
\end{proof}

Putting matters together this concludes the proof of Theorem \ref{t=A}.

\end{document}